\documentclass[journal]{IEEEtran}
\IEEEoverridecommandlockouts 

\usepackage{cite}
\usepackage{amsmath, amstext, amssymb}
\usepackage{array}
\usepackage{stfloats}
\usepackage{url}
\usepackage{graphicx}

\usepackage[]{hyperref}

\begin{document}
\bstctlcite{IEEEexample:BSTcontrol}
\renewcommand{\baselinestretch}{0.96}
%
\title{Chance-Constrained AC Optimal Power Flow: \\ Reformulations and Efficient Algorithms}

\author{Line~Roald~
        and~G\"oran~Andersson~


%

\thanks{}
\thanks{L. Roald is with CNLS and T Division of Los Alamos National Laboratory, Los Alamos, NM, United States. G. Andersson is with the Power System Laboratory, ETH Zurich, Switzerland. E-mail: roald@lanl.gov}}



\maketitle

\begin{abstract}
Higher levels of renewable electricity generation increase uncertainty in power system operation. To ensure secure system operation, new tools that account for this uncertainty are required. In this paper, we adopt a chance-constrained AC optimal power flow formulation, which guarantees that generation, power flows and voltages remain within their bounds with a pre-defined probability. We then discuss different chance-constraint reformulations and solution approaches for the problem. 
We first describe an analytical reformulation based on partial linearization, which enables us to obtain a tractable representation of the optimization problem. 
We then provide an efficient algorithm based on an iterative solution scheme which alternates between solving a deterministic AC OPF problem and assessing the impact of uncertainty. The flexibility of the iterative scheme enables not only scalable implementations, but also alternative chance-constraint reformulations. In particular, we suggest two sample based reformulations that do not require any approximation or relaxation of the AC power flow equations. 

In a case study based on four different IEEE systems, we assess the performance of the method, and demonstrate scalability of the iterative scheme. We further show that the analytical reformulation accurately and efficiently enforces chance constraints in both in- and out-of-sample tests, and that the analytical reformulations outperforms the two alternative, sample-based chance constraint reformulations.

\end{abstract}


\section{Introduction}
Over the last decade, 
energy production from wind and solar power 
has reached significant levels in many countries across the world. 
For power system operators, the integration of renewable energy poses a variety of challenges, from long-term generation adequacy to a reduction of system rotational inertia. 
In this paper, we address the question of how to assess and mitigate the impact of forecast uncertainty from renewable generation in day-to-day operational planning. 

A fundamental tool in the operational planning is the Optimal Power Flow (OPF), an optimization problem commonly involved in market clearing and security assessment processes \cite{stott2012whitepaper, capitanescu2011stateoftheart}.  
Most OPF problems aim at minimizing operational cost while ensuring secure operation by enforcing constraints such as transmission capacity, voltage and generation limits.
Traditionally, the OPF has been formulated as a \emph{deterministic} problem.
However, with more renewable generation, it becomes increasingly important to account for forecast uncertainty and treat the OPF as a \emph{stochastic} problem. 

Consequently, a wide variety of approaches and methods to account for uncertainty within the OPF have recently been proposed in literature. These include, among others, robust and worst-case methods \cite{panciatici2010, capitanescu2012, warrington2013, lorca2015}, two- and multi-stage stochastic programming based on samples \cite{bouffard2008, morales2009, papavasiliou2013, mueller2014} or stochastic approximation techniques \cite{tinoco2012}, and chance-constrained formulations 
\cite{vrakopoulou2013, roald2013, bienstock2014, zhang2011, li2011, hojjat2015, baker2016, vrakopoulou2013AC, schmidli2016, venzke17arxiv}.
In this paper, we work with chance constraints, which ensure that the system constraints will be satisfied with a specified probability. 
Discussions with transmission system operators \cite{umbrella} have revealed that choosing an acceptable violation probability is perceived as an intuitive and transparent way of determining a probabilistic security level. 
Chance constraints also align well with several industry practices, such as the probabilistic reserve dimensioning in ENTSO-E \cite{ENTSOE2012, ENTSOE2013supportingLFCR} 
or the definition of reliability margins in the European market coupling \cite{CWE2011}. 

While chance constraints offer an intuitive way of limiting risk from forecast uncertainty, the resulting optimization problems are generally hard to solve. 
Most literature has so far considered the linear DC power flow approximation, since modelling the impact of uncertainty in the full, non-linear AC power flow problem is a challenge. However, many applications, such as distribution system optimization or transmission system security assessment, require the more accurate AC power flow equations. Modelling AC power flow also allows us to include and probabilistically enforce a range of new constraints, such as voltage limits, reactive power constraints and transmission constraints based on apparent power or current magnitude. 
Formulations that have attempted to solve the full Chance-Constrained AC OPF (AC CC-OPF) include \cite{zhang2011, li2011, hojjat2015, baker2016, vrakopoulou2013AC, schmidli2016}. 
In \cite{zhang2011}, the problem is solved using an iterative gradient calculation and numerical integration, while \cite{li2011} employs an iterative approach based on the cumulant method and Cornish-Fischer expansion. 
AC CC-OPF based on linearized equations has been used to find optimal redispatch schedules \cite{hojjat2015} and solve voltage-constrained OPF in distribution grids \cite{baker2016}. 
Modelling of the AC CC-OPF based on convex relaxations of the AC power flow equations have been applied in \cite{vrakopoulou2013AC, venzke17arxiv}.
The probably most comprehensive AC CC-OPF formulation to date is provided in \cite{vrakopoulou2013AC}, which is based on the SDP relaxation for the power flows and a sample-based approach to reformulate the chance constraint. 
However, checking that the relaxed power flow solution is physically meaningful (i.e., that the relaxation is exact) for all points within the uncertainty set is still a challenge that need to be addressed. Furthermore, the high computational requirements limits the approach to small systems.
Indeed, all of the above AC CC-OPF methods have only been applied to small test cases, 
signalling a need for scalable approaches \cite{capitanescu2016}.

In this paper, we discuss different reformulations and solution algorithms for the AC CC-OPF, with the aim of providing both an accurate representation of the AC power flow constraints and maintain scalability to large instances. 
We first introduce a detailed model of the AC CC-OPF with separate chance constraints, based on a classical model for power system operation where system balance is maintained using Automatic Generation Control (AGC) and reactive power is adjusted using local voltage regulation at PV buses. The model is an extension of the model used in \cite{schmidli2016} which neglected the impact of uncertainty on generation constraints. 

Second, we present an analytical reformulation based on a \emph{partial} linearization of the AC power flow equations. 
This reformulation includes the full, non-convex AC power flow equations for the forecasted operating point, without any approximation or relaxation, and hence guarantees AC feasibility for this operating point. Assuming that the forecast errors are small, the impact of uncertainty is modelled using a linearization around this point. 
Note that this approach is more accurate than the full linearization in \cite{zhang2011, li2011, hojjat2015, baker2016}, and does not require a relaxation as in \cite{vrakopoulou2013AC}. Similar partial linearization methods were applied to stochastic load flow, see \cite{dopazo1975} and discussion therein, and were also adopted in the risk-based OPF in \cite{fu2001}. 
The partial linearization approximates the impact of uncertainty as a linear function of the uncertain power injections. This linear relationship enables the application of the analytical chance constraint reformulation first proposed for DC power flow in \cite{roald2013, bienstock2014}, which is applicable to both Gaussian and  partially unknown uncertainty distributions \cite{roaldArxiv}. The analytical chance-constraint reformulation was developed for apparent power constraints in \cite{qu2015}, and preliminary results for the application to the AC CC-OPF were presented in \cite{schmidli2016}. Compared to previous work, this paper extends the power system model to include chance-constraints on generators, and provides a more in-depth analysis of the performance of the analytical reformulation. In particular, we assess the performance of the Gaussian reformulation in systems with non-Gaussian uncertainty and demonstrate the benefit of optimizing the system response to uncertainty. 

The analytical reformulation leads to closed form expressions for the constraints and can hence be solved for as a one-shot optimization problem. Due to the closed-form expressions, the number of constraints with the analytical reformulation remain the same as in the deterministic case, as opposed to scenario based methods such as \cite{vrakopoulou2013AC}, where the required number of samples significantly increases with the size of the problem. However, the reformulated constraints are complex, non-convex functions of the decision variables, and therefore still might be computationally prohibitive. As an alternative solution algorithm, we use the iterative solution scheme implemented on a small test case in \cite{schmidli2016}. This iterative scheme is based on the observation that the impact of uncertainty can be interpreted as a tightening of the original, deterministic constraints. The algorithm hence alternates between (i) solving a deterministic AC OPF with a fixed constraint tightening and (ii) assessing the impact of uncertainty and computing the required tightening for the obtained optimal solution. The algorithm is deemed to have converged when the tightenings remain the same between subsequent iterations. 

In this paper, we first show that the iterative approach provides a solution which is similar to the local optimum of the one-shot optimization, verifying that the iterative approach works well. Second, we demonstrate scalability of the iterative approach, which utilizes the efficiency of existing deterministic AC OPF solvers, by solving the large-scale Polish test case with 941 uncertain loads in 32s on a desktop computer.

By decoupling the uncertainty assessment from the solution of the AC OPF, we observe that the iterative algorithm not only admits scalable implementations, but also enable more general and accurate reformulations of the chance constraints without compromising computational tractability. Any method can be used to assess the uncertainty impact, including methods that would be computationally intractable if included directly in the optimization problem. 
In this paper, we suggest two alternative sample-based reformulations based on (i) a Monte-Carlo simulation and (ii) the scenario approach \cite{calafiore2006}. The advantage of the sample-based methods is that they neither require any assumption about the distribution (except for availability of samples), nor approximation or relaxation of the AC power flow equations. Both methods would be a computational challenge to include directly in the problem, but are easily incorporated using the iterative framework, as demonstrated in the case study. We use the implementation of the sample-based reformulations to benchmark the performance of the analytical reformulation in the case study.

In summary, the main contributions of this paper are the following:
\begin{enumerate}
    \item We extend the modelling of the AC CC-OPF in \cite{schmidli2016} to include chance constraints on generator active and reactive power outputs, and optimizing the system response to uncertainty.
    \item We compare the solutions and solution times obtained when solving the AC CC-OPF either as a one-shot optimization problem or using an iterative solution algorithm, and show that the iterative algorithm provides optimal results in the considered cases. 
    \item We demonstrate the scalability of the iterative algorithm by solving the Polish test case with 941 uncertain loads in 32s on a desktop computer.
    \item We show that the iterative algorithm can be used to include more general and accurate chance-constraint reformulations without compromising computational tractability, by proposing and implementing two sample-based chance constraint reformulations based on a Monte-Carlo simulation and the scenario approach.
    \item We demonstrate that the analytical reformulation based on a Gaussian uncertainty is accurate both for in- and out-of-sample tests, and outperforms the sample-based reformulations.
\end{enumerate}

The remainder of the paper is structured as follows. We describe the uncertainty modelling in Section \ref{sec:Uncertainty}, and the AC OPF formulation with chance constraints in Section \ref{sec:ACCCOPF}. The analytical chance constraint reformulation is described in Section \ref{sec:Analytical}. Section \ref{sec:Solution} describes the one-shot and iterative solution approaches, while Section \ref{sec:Alternative} discusses alternative chance constraint reformulations. In Section \ref{sec:Case}, we present the case study to assess the performance of the proposed method, before Section \ref{sec:Practical} discusses aspects relevant to the practical implementation of AC CC-OPF in industry. Section \ref{sec:Conclusion} summarizes and concludes.

\section{Power System Modelling under Uncertainty}
\label{sec:Uncertainty}
We consider a power system where $\mathcal{N},~\mathcal{L}$ denote the set of nodes and lines, and  $|\mathcal{N}| = m$ and $|\mathcal{L}| = l$. 
The set of nodes with uncertain demand or production of energy is given by $\mathcal{U}\subseteq\mathcal{N}$, while the set of conventional generators is denoted by $\mathcal{G}$. 
To simplify notation, we assume that there is one conventional generator with active and reactive power outputs $p_{G,i},~q_{G,i}$, one composite uncertainty source  $p_{U,i},~q_{U,i}$ and one demand $p_{D,i},~q_{D,i}$ per node, such that $|\mathcal{G}|=|\mathcal{U}|=|\mathcal{N}| = m$.
Nodes without generation or load can be handled by setting the respective entries to zero, and nodes with multiple entries can be handled through a summation. 
To model generation and control at different types of buses, we add subscripts $_{PQ},~_{PV},$ and $_{\theta V}$ to distinguish between PQ, PV and $\theta V$ (reference) buses. 
In our formulation, we assume only one $\theta V$ bus, although the method could be extended to account for a distributed slack bus.

\subsubsection{Uncertainty Modelling}
The deviations in demand or production at any given node can be due to, e.g., load fluctuations, forecast errors for wind or PV or intra-day electricity trading. 
We model the deviations in active power as the sum of the forecasted value $p_U$ and a zero mean fluctuation $\omega$,
\begin{eqnarray}
	\tilde{p}_{U} = p_U + \omega~. \nonumber
\end{eqnarray}
We assume a \emph{constant power factor} $\cos \phi$ for the uncertain injections, such that the reactive power injections are given by
\begin{eqnarray}
	\tilde{q}_U = q_U + \gamma\omega~, \quad \text{where } 
    \gamma = \sqrt{\left(1-\cos^2 \phi\right)/\left(\cos^2 \phi\right)}~. \nonumber
\end{eqnarray}
The variable $\gamma$ will be referred to as the power ratio in the following. 
The reactive power injections could also be modelled in other ways, e.g. assuming that the reactive power injections remain constant, that the reactive power can be dispatched (at least partially) independent of the active power production or that some uncertainty sources (e.g., large wind farms) participate in controlling the voltage at their point of connection. These types of control can be included in the formulation without any conceptual changes. 

\subsubsection{Generation and Voltage Control}
The generation dispatch $p_G,~q_G$ and voltage magnitudes $v$ are scheduled by the system operator. The forecasted operating point, corresponding to the situation with $\omega=0$, is assumed to be a balanced operating point which satisfies the nodal power balance equations. 
The controllable generators further adjust their reactive and active power outputs to ensure power balance and maintain the desired voltage profile during fluctuations $\omega\neq 0$. 
We assume that active power generation is balanced by the Automatic Generation Control (AGC) \cite{borkowska}, and that the power generation mismatch $\Omega=\sum_{i\in\mathcal{U}} \omega_i$ is divided among generators according to participation factors $\alpha$. 
While the deviation $\Omega$ due to uncertainty is the main source of power imbalance in the system, there is also an additional power mismatch due to changes in the active power losses. The total change in active power loss, which is a-priori unknown and denoted by $\delta p$, is assumed to be balanced by the generator at the reference bus\footnote{The balancing of the active power losses could also be done directly through the AGC by redefining $\hat\Omega=\sum_{i\in\mathcal{U}} \omega_i + \delta p$. However, since the change in losses $\delta p$ is a secondary effect (arising because power flow changes due to the power deviations $\omega$) and was found to be small in simulations, we assume balancing by the generator at the reference bus only.}. 
\begin{subequations}
\label{eq:generatoradjustment}
\begin{align}
	&\tilde{p}_{G,i}(\omega) = p_{G,i} - \alpha_i \Omega, \quad &&\forall_{i\in\mathcal{G}_{PQ},~\mathcal{G}_{PV}}, \\
	&\tilde{p}_{G,i}(\omega) = p_{G,i} - \alpha_i \Omega + \delta p, \quad &&\forall_{i\in\mathcal{G}_{\theta V}}, \\
	&\sum_{i\in\mathcal{G}}\alpha_i = 1. 
\end{align}
\end{subequations}
We further assume that reactive power is controlled locally at PV and $\theta$V buses (where generators change their outputs by $\delta q$ to keep voltage magnitudes constant), whereas generators at PQ buses keep their output constant,
\begin{align}
	\tilde{q}_{G,i}(\omega) &= q_{G,i} + \delta q_i, \quad &&\forall_{i\in\mathcal{G}_{PV},~\mathcal{G}_{\theta V}}, \nonumber \\
	\tilde{q}_{G,i}(\omega) &= q_{G,i}, \quad &&\forall_{i\in\mathcal{G}_{PQ}}. \nonumber
\end{align}
Conversely, the voltage magnitude is fixed at the reference and PV buses, but varies at PQ buses:
\begin{align}
	\tilde{v}_j(\omega) &= v_{j},            \quad &&\forall_{j\in\mathcal{N}_{PV},~{N}_{\theta V}},\nonumber \\
	\tilde{v}_j(\omega) &= v_{j} + \delta v_j, \quad &&\forall_{j\in\mathcal{N}_{PQ}}. \nonumber
\end{align}
Here, $v_j$ represents the voltage magnitude at the $j^{th}$ node.
Note that a centralized Automatic Voltage Regulation (AVR) scheme similar to the one proposed in \cite{vrakopoulou2013AC}, with control variables analogous to the participation factors $\alpha$ of the AGC, could also be incorporated.

\subsubsection{Power flows}
With deviations in the active and reactive power injections, the power flows throughout the system also change. 
Assuming a thermally constrained power system, we state transmission constraints in terms of current constraints, with $i_{ij}$ denoting the magnitude of the current flow from bus $i$ to bus $j$. Changes in the power injections induce a change $\delta i_{ij}$ in the current magnitudes, such that the total current magnitude $\tilde{i}_{ij}$ is given by 
\begin{align}
	\tilde{i}_{ij}(\omega) &= i_{ij} + \delta i_{ij}, \quad &&\forall_{ij\in\mathcal{L}}.
\end{align}

\section{Chance-Constrained AC Optimal Power Flow}
\label{sec:ACCCOPF}
We formulate the full, original AC CC-OPF as 
\begin{subequations}
\label{originalACCCOPF}
\begin{align}
\min_{p_G, q_G, v, \theta, \gamma, \alpha} ~~&  \sum_{i\in \mathcal{G}}\left(c_{2,i}p_{G,i}^2 + c_{1,i}p_{G,i} +c_{0,i}\right)&& \label{ac_obj}\\
\text{s.t.}  ~~         
&f\left( \tilde{\theta}(\omega), \tilde{v}(\omega), \tilde{p}(\omega),\tilde{q}(\omega) \right) = 0 ~&&\forall_{\omega \in \mathcal{D}} \label{ac_powerbal}\\
&\mathbb{P}(\tilde{p}_{G,i}(\omega) \leq p_{G,i}^{max})\geq 1-\epsilon_P, ~~&&\forall_{i\in\mathcal{G}} \label{ac_pMax}\\
&\mathbb{P}(\tilde{p}_{G,i}(\omega) \geq p_{G,i}^{min})\geq 1-\epsilon_P, ~~&&\forall_{i\in\mathcal{G}} \label{ac_pMin}\\
&\mathbb{P}(\tilde{q}_{G,i}(\omega) \leq q_{G,i}^{max})\geq 1-\epsilon_Q, ~~&&\forall_{i\in\mathcal{G}} \label{ac_qMax}\\
&\mathbb{P}(\tilde{q}_{G,i}(\omega) \geq q_{G,i}^{min})\geq 1-\epsilon_Q, ~~&&\forall_{i\in\mathcal{G}} \label{ac_qMin}\\
&\mathbb{P}(\tilde{v}_{j}(\omega) \leq v_{j}^{max})\geq 1-\epsilon_V, ~~&&\forall_{j\in\mathcal{N}} \label{ac_vMax}\\
&\mathbb{P}(\tilde{v}_{j}(\omega) \geq v_{j}^{min})\geq 1-\epsilon_V, ~~&&\forall_{j\in\mathcal{N}} \label{ac_vMin}\\
&\mathbb{P}(\tilde{i}_{ij}(\omega) \leq i_{ij}^{max})\geq 1-\epsilon_I, ~~&&\forall_{ij\in\mathcal{L}} \label{ac_iMax}\\
&\theta_{\theta V} = 0 \label{ac_ref} 
\end{align}
\end{subequations}
The objective \eqref{ac_obj} minimizes the cost of active power generation, with $c_2$, $c_1$ and $c_0$ being the quadratic, linear and constant cost coefficients. 
Eq. \eqref{ac_powerbal} represent the AC power balance constraints for all possible $\omega$ within the uncertainty set $\mathcal{D}$. These are functions of the nodal voltage angles $\tilde{\theta}(\omega)$ and magnitudes $\tilde{v}(\omega)$, as well as the nodal injections of active $\tilde{p}(\omega)$ and reactive power $\tilde{q}(\omega)$, and are given by
\begin{subequations}
\label{pf_polar}
\begin{align}
\nonumber
\tilde{p}_j(\omega) = \tilde{v}_j(\omega) \sum_{k=1}^n \tilde{v}_k(\omega&)  \left(\mathbf{G}_{jk}\cos\left(\tilde{\theta}_j(\omega) - \tilde{\theta}_k(\omega)\right) \right. \\ 
& \left. +\, \mathbf{B}_{jk}\sin\left(\tilde{\theta}_j(\omega) - \tilde{\theta}_k(\omega)\right) \right) \\
\nonumber
\tilde{q}_j(\omega) = \tilde{v}_j(\omega) \sum_{k=1}^n \tilde{v}_k(\omega&)   \left(\mathbf{G}_{jk}\sin\left(\tilde{\theta}_j(\omega) - \tilde{\theta}_k(\omega)\right) \right. \\ & \left. -\, \mathbf{B}_{jk}\cos\left(\tilde{\theta}_j(\omega) - \tilde{\theta}_k(\omega)\right) \right)
\end{align}
\end{subequations}
where $\mathbf{G}$ and $\mathbf{B}$ denote the real and imaginary components, respectively, of the network admittance matrix.
The nodal power injections are the sum of the power injections from generator, loads and uncertainty sources,
\begin{align}
&\tilde{p}(\omega) = \tilde{p}_G(\omega) + p_D + p_U + \omega, \\
&\tilde{q}(\omega) = \tilde{q}_G(\omega) + q_D + q_U + \gamma\omega. 
\end{align}

The remaining constraints are generation constraints for active and reactive power \eqref{ac_pMax}-\eqref{ac_qMin}, constraints on the voltage magnitudes at each bus \eqref{ac_vMax}, \eqref{ac_vMin} and transmission constraints on the current magnitudes \eqref{ac_iMax}. All these constraints are formulated as chance constraints with acceptable violation probabilities of $\epsilon_P,~\epsilon_Q, ~\epsilon_V$ and $\epsilon_I$, respectively. 
Finally, the voltage angle at the reference bus is set to zero by \eqref{ac_ref}.

\subsubsection{Interpretation of the Chance Constraints}
The chance constraints \eqref{ac_pMax}-\eqref{ac_iMax} limit the probability of constraint violation, but do not describe what happens in cases where the constraints are violated. 
For an interpretation of the meaning of violations, we can distinguish between two main types of constraints: 
\begin{enumerate}
    \item \emph{Hard} constraints such as generation constraints \eqref{ac_pMax}-\eqref{ac_qMin}, which are physically impossible to violate. Violations of hard generation constraints would indicate a situation in which a generator reaches an upper or lower limit and is unable to provide reactive power support or balancing power according to the prescribed participation factor $\alpha_i$. In this case, the automatic control actions such as the AGC might leave the system imbalanced and manual intervention would be needed\footnote{Note that the generator adjustment to the uncertainty \eqref{eq:generatoradjustment} only depend on the total power imbalance $\Omega$ (except for the slack bus, which has an additional component balancing the change in losses). The total imbalance is a scalar random variable, which means that the generator constraints required to guarantee availability of AGC with a probability of $\epsilon_P$ as in Eqs. \eqref{ac_pMax}, \eqref{ac_pMin} are given by
    \begin{equation}
    p-\alpha \Omega_{1-\epsilon_P} \leq p^{max} \nonumber
    \end{equation}
    where $\Omega_{1-\epsilon_P}$ represents the $(1-\epsilon_P)$ quantile of the total imbalance. This implies that if $|\Omega| \leq \Omega_{1-\epsilon_P}$, all generator constraints will be \emph{jointly} satisfied (i.e., none of the generators will reach their limit). The violation probability $\epsilon_P$ thus represent the probability of insufficient AGC capacity, and the generator chance constraints can be interpreted as a probabilistic reserve criteria, similar to criteria enforced in industry\cite{roald2016corrective}. As an example, Swissgrid (the swiss system operator) enforces $\epsilon_P\leq 0.001$ in their reserve procurement process \cite{abbaspourtorbati2016}.}. 
    \item \emph{Soft} constraints such as voltage and current magnitude constraints \eqref{ac_vMax}-\eqref{ac_iMax}, where a constraint violation implies an actual under- or over-voltage, or a transmission line overload. Constraint violations of soft constraints might either be tolerated if the magnitude and duration are not too large, or removed through additional control actions such as generation redispatch.
\end{enumerate}
In general, we can interpret the chance constraint as the probability with which the operator needs to implement additional control actions to secure the system in real time. Choosing a low acceptable violation probability will make system operation safer and less stressful for the operator, but also more costly. Choosing a high acceptable violation probability is risky, as it might require frequent application of real-time control and there is no guarantee that such controls will be available. 


\subsubsection{Optimization of System Response} The AC CC-OPF formulation includes the AGC participation factors $\alpha$ and the power ratio $\gamma$ as optimization variables, and hence assumes that these are controllable. In many systems, the values for $\alpha,~\gamma$ are either not controllable or are pre-determined based on, e.g., separate reserve procurement processes or generator physical constraints. In this case, the variables $\alpha,~\gamma$ are reduced to pre-specified problem parameters. In this paper, we will mostly assume those variables are pre-specified, but will also assess the value of being able to control them.

\section{Analytical Reformulation of the Chance Constraints}
\label{sec:Analytical}
The problem \eqref{originalACCCOPF} is not tractable in its current form. First, \eqref{ac_powerbal} is semi-infinite, as the set of fluctuations $\mathcal{D}$ is uncountable. Second, the chance constraints \eqref{ac_pMax}-\eqref{ac_iMax} must be reformulated into tractable constraints that can be handled by an optimization solver.
To obtain a tractable optimization problem, we suggest a reformulation of \eqref{ac_obj}-\eqref{ac_ref} based on the following main ideas: 
\begin{itemize}
	\item For the forecasted operating point, given by $\omega = 0$ and the scheduled nodal power injections $p,~q$, we solve the full AC power flow equations. This ensures an accurate, AC feasible solution for the forecasted system state.   
	\item The impact of the uncertainty $\omega$ is modelled using a first-order Taylor expansion around the forecasted operating point. This approximation can be expected to be accurate when the forecast errors $\omega$ are small, which is a reasonable assumption close to real time operation when forecasts tend to be good\footnote{Note that we are not assuming that the renewable generation itself is small, but rather that the \emph{forecast error} is small. Therefore, we expect the method to provide accurate results in any system with good forecasting procedures, regardless of the renewable penetration.}.  
\end{itemize}
Following the above ideas, we replace \eqref{ac_powerbal} by a single set of deterministic equations for the forecasted operating point,
\begin{equation}
f\left( \theta, v, p, q \right) = 0, \label{nominal_OP}
\end{equation}
where $\theta,~v$ are the voltage angles and magnitudes corresponding to the scheduled injections $p,~q$. In the following, we will denote this operating point by $\mathbf{x}=\left( \theta, v, p, q \right)$. 
We then define sensitivity factors with respect to the fluctuations $\omega$, i.e.,
\begin{align}
&\mathbf{\Gamma}_{P}(\mathbf{x},\!\alpha,\!\gamma) \!\!= \!\!\!\!\left.\begin{array}{c}\frac{\partial p_G}{\partial \omega}\end{array} \!\!\! \right\vert_{(\omega=0, \mathbf{x}, \alpha, \gamma)},\!
&&\mathbf{\Gamma}_{Q}(\mathbf{x},\!\alpha,\!\gamma) \!\!= \!\!\!\!\left.\begin{array}{c}\frac{\partial q_G}{\partial \omega}\end{array} \!\!\! \right\vert_{(\omega=0, \mathbf{x}, \alpha, \gamma)}\nonumber\\
&\mathbf{\Gamma}_{V}(\mathbf{x},\!\alpha,\!\gamma) \!\!= \!\!\!\!\left.\begin{array}{c}\frac{\partial v}{\partial \omega} \end{array} \!\!\! \right\vert_{(\omega=0, \mathbf{x}, \alpha, \gamma)},\!
&&\mathbf{\Gamma}_{I}(\mathbf{x},\!\alpha,\!\gamma) \!\!= \!\!\!\!\left.\begin{array}{c}\frac{\partial i}{\partial \omega} \end{array} \!\!\! \right\vert_{(\omega=0, \mathbf{x}, \alpha, \gamma)} 
\nonumber
\end{align}
Here, $\mathbf{\Gamma}_P,~\mathbf{\Gamma}_Q$ denote the sensitivity factors for the active and reactive power output of the conventional generators, while $\mathbf{\Gamma}_V,~\mathbf{\Gamma}_I$ are the sensitivity factors for the voltage and current magnitudes. While the general approach of deriving sensitivity factors is well known in literature, see e.g. \cite{dopazo1975}, 
a detailed derivation of  $\mathbf{\Gamma}_P,~\mathbf{\Gamma}_Q,~\mathbf{\Gamma}_V$ and $\mathbf{\Gamma}_I$, including considerations related to the handling of PV, PQ and $\theta$V buses, can be found in \cite{roald2016thesis}. Note that the sensitivity factors depend non-linearly on the forecasted operating point $\mathbf{x}$, and linearly on the participation factors $\alpha$ and the power ratio $\gamma$.

The sensitivity factors allow us to approximate the chance constraints 
as linear functions of  the random variables $\omega$,
e.g., for the current constraint \eqref{ac_iMax} we obtain
\begin{equation}
	\mathbb{P}(i_{ij} + \mathbf{\Gamma}_{I(ij,\cdot)} \omega \leq i_{ij}^{max})\geq 1-\epsilon_I, ~~\forall_{ij\in\mathcal{L}} \label{ac_iMax_lin}
\end{equation}
where $\mathbf{\Gamma}_{I(ij,\cdot)}$ is a row vector representing the sensitivity of the current on line $ij$ to changes in $\omega$. 
The linear dependence on $\omega$ enables the use of an analytical chance constraint reformulation \cite{roald2013, qu2015, schmidli2016}, 
even though the sensitivity factors $\mathbf{\Gamma}_P,~\mathbf{\Gamma}_Q,~\mathbf{\Gamma}_V,~\mathbf{\Gamma}_I$ and the current magnitudes $i_{ij}$ are non-linear functions of the decision variables. 
Assuming that the fluctuations $\omega$ follows a multivariate normal distribution with zero mean and covariance matrix $\Sigma_W$, we obtain the following expression for \eqref{ac_iMax_lin} \cite{qu2015},
\begin{equation}
	i_{ij} + \Phi^{-1}(1\!-\!\epsilon_I) \!\parallel\! \mathbf{\Gamma}_{I(ij,\cdot)}\Sigma_W^{1/2} \!\parallel_2 \leq i_{ij}^{max}. 
	\label{curr_uncmarg}
\end{equation}
Here, $\Phi^{-1}(1-\epsilon_I)$ represents the inverse cumulative distribution function of the standard normal distribution, evaluated at $1-\epsilon_I$. Note that the method can guarantee safety for more general distributions by replacing the inverse cumulative distribution function of the normal distrbution by general probabilistic bounds, which only require knowledge of the mean and covariance of the random variables \cite{roaldArxiv}.

We observe that the consideration of uncertainty introduce an \emph{uncertainty margin}, i.e., a tightening of the constraint on the forecasted current which is necessary to secure the system against uncertainty. Denoting this uncertainty margin by $\lambda_I$, we rewrite \eqref{curr_uncmarg} as
\begin{align}
    &i_{ij} \leq i_{ij}^{max} -\lambda_I(\mathbf{x},\alpha,\gamma), \quad \text{with}
	\label{currentII}\\
	&\lambda_I(\mathbf{x},\alpha,\gamma) = \Phi^{-1}(1\!-\!\epsilon_I) \!\parallel\! \mathbf{\Gamma}_{I(ij,\cdot)}\Sigma_W^{1/2} \!\parallel_2 
	\label{currentIII}
\end{align}
The uncertainty margin $\lambda_I$ depends on non-linearly on both the forecasted operating point $\mathbf{x}$, the participation factors $\alpha$ and the power ratio $\gamma$. 

\subsection{Reformulated Chance-Constrained Problem}
Applying the analytical reformulation to all chance constraints and using the definition of the uncertainty margins from above, we express the AC CC-OPF \eqref{originalACCCOPF} as
\begin{subequations}
\label{approxACCCOPF}
\begin{align}
\min_{\mathbf{x}, p_G, q_G, \alpha, \gamma} ~&  \sum_{i\in \mathcal{G}}\left(c_{2,i}p_{G,i}^2 + c_{1,i}p_{G,i} +c_{0,i}\right)&  \label{eq1} \\
\text{s.t.} ~~          
&f\left( \theta, v, p, q \right) = 0   \label{eq2}\\
& p_G^{min} + \lambda_P \leq p_G \leq p_G^{max} - \lambda_P \\
& q_G^{min} + \lambda_Q \leq q_G \leq q_G^{max} - \lambda_Q \\
& v^{min} + \lambda_V \leq v \leq v^{max} - \lambda_V \\
& i \leq i^{max} - \lambda_I \\
&\theta_{\theta V} = 0 \label{eq14}
\end{align}
\end{subequations}
The uncertainty margins for currents $\lambda_I$ are defined by \eqref{currentIII}, 
while the voltage uncertainty margins $\lambda_V$ are given by
\begin{subequations}
\begin{align}
& \lambda_{V,j} = \Phi^{-1}(1-\epsilon_V) \!\parallel\! \mathbf{\Gamma}_{V(j,\cdot)}\Sigma_W^{1/2}  \!\parallel_2
&&\forall_{j\in\mathcal{N}_{PQ}} \label{uncV}, \\
& \lambda_{V,j} = 0 
&&\forall_{j\in\mathcal{N}_{PV}, \mathcal{N}_{\theta V}} \label{uncV}. 
\end{align}
and the $\lambda_P,~\lambda_Q$ for active and reactive power are given by
\label{lambdas}
\begin{align}
& \lambda_{P,i} = \alpha_i \Phi^{-1}(1-\epsilon_P)\sigma_\Omega,
&&\!\!\!\!\!\!\!\!\!\!\!\!\!\!\!\!\!\!\!\!\!\!\!\!\!\!\!\!\!\!\!\!\!\!\forall_{i\in\mathcal{G}_{PQ},~ \mathcal{G}_{PV}}, \label{uncP}\\
& \lambda_{P,i} =  \Phi^{-1}(1-\epsilon_P) \!\parallel\! (-\alpha_i \mathbf{1}_{1,m} + \mathbf{\Gamma}_{P(i,\cdot)})\Sigma_W^{1/2}  \!\parallel_2,&&\\ & &&\!\!\!\!\!\!\!\!\!\!\!\!\!\!\!\!\!\!\!\!\!\!\!\!\!\!\!\!\!\!\!\!\!\!\forall_{i\in\mathcal{G}_{\theta V}} \nonumber\\
& \lambda_{Q,i} = \Phi^{-1}(1-\epsilon_Q) \!\parallel\! \mathbf{\Gamma}_{Q(i,\cdot)}\Sigma_W^{1/2}  \!\parallel_2,
&&\!\!\!\!\!\!\!\!\!\!\!\!\!\!\!\!\!\!\!\!\!\!\!\!\!\!\!\!\!\!\!\!\!\!\forall_{i\in\mathcal{G}_{PV},~ \mathcal{G}_{\theta V} }, \\
& \lambda_{Q,i} = 0,
&&\!\!\!\!\!\!\!\!\!\!\!\!\!\!\!\!\!\!\!\!\!\!\!\!\!\!\!\!\!\!\!\!\!\!\forall_{i\in\mathcal{G}_{PQ}}.
\end{align}

\end{subequations}

Here, $\sigma_\Omega$ represents the standard deviation of the total active power imbalance $\Omega$, given by $\sigma_{\Omega}=\parallel \! 1^T \Sigma^{1/2} \! \parallel_2$. The vector $\mathbf{1}_{1,m}$ is a row vector with $m$ elements and entries equal to 1.

\section{Solution Algorithms}
\label{sec:Solution}
We now discuss two solution algorithms for the AC CC-OPF described above. 

\subsection{One-Shot Optimization} 
The reformulated AC CC-OPF given by \eqref{approxACCCOPF}, \eqref{currentIII}, \eqref{lambdas} is a continuous, non-convex optimization problem, which can be solved directly by a suitable non-linear solver. 
The solver optimizes the scheduled generation dispatch $p_G,~q_G$, while inherently accounting for the dependency of the sensitivity factors $\mathbf{\Gamma}_P,~\mathbf{\Gamma}_Q,~\mathbf{\Gamma}_V$ and $\mathbf{\Gamma}_I$ on the forecasted operating point $\mathbf{x}$, the participation factors $\alpha$ and the power ratios $\gamma$. 
By including $\alpha$ and $\gamma$ as optimization variables, it is thus possible to optimize not only the scheduled dispatch, but also the procurement of reserves and voltage control during deviations. Note that the problem is not guaranteed to converge to a global optimum, as it is non-convex.

The drawback of attempting a one-shot solution is the problem complexity, which might lead to long solution times.
The deterministic AC OPF is already a non-convex problem, and adding additional terms with complex dependencies on the decision variables only increases the difficulty of obtaining a solution. This can be a bottleneck for adoption in more realistic settings, where scalability and robustness of the OPF are important criteria \cite{stott2012whitepaper}. 

\subsection{Iterative Solution Algorithm}
\label{sec:iterative}
To address the problem of increased computational complexity, we apply the iterative solution algorithm from \cite{schmidli2016}, which allows us to obtain a solution to the AC CC-OPF given by \eqref{approxACCCOPF}, \eqref{currentIII}, \eqref{lambdas} using any existing AC OPF tool. 
The iterative solution algorithm is based on the observation that the impact of uncertainty is only visible through the uncertainty margins $\lambda$. The scheme thus alternates between solving a deterministic AC OPF given by \eqref{approxACCCOPF} with fixed $\lambda^{\kappa-1}$ to obtain $\mathbf{x}^\kappa$, and evaluating $\lambda^{\kappa}=\lambda(\mathbf{x}^\kappa,\alpha,\gamma)$ based on \eqref{currentIII}, \eqref{lambdas}. When the uncertainty margins $\lambda^\kappa$ do not change between iterations, the algorithm has converged to a feasible solution. 
%
\\
More specifically, the algorithm consist of the following steps:
\vspace{+6pt}
\hrule
\vspace{+2pt}
\begin{enumerate}
	\item \emph{Initialization:} Set uncertainty margins $\lambda_P^0=\lambda_Q^0=\lambda_V^0=\lambda_I^0=0$, and iteration count $\kappa=1$.
	\vspace{+3pt}
	\item \emph{Solve AC OPF:} Solve the deterministic AC OPF defined by \eqref{approxACCCOPF} with fixed $\lambda^{\kappa-1}$, and obtain a solution for the forecasted operating point $\mathbf{x}^\kappa$. 
	\vspace{+3pt}
	\item \emph{Evaluate uncertainty margins:} Compute the uncertainty margins of the current iteration $\lambda_P^{\kappa},~\lambda_Q^{\kappa},~\lambda_V^{\kappa}$ and $\lambda_I^{\kappa}$ based on $\mathbf{x}^\kappa$. Then, evaluate the maximum change in the uncertainty margins from the last iteration,
	\begin{align}
	& \eta^{\kappa}_P = \max\{|\lambda_P^{\kappa}-\lambda_P^{\kappa-1}|\}, \quad
	\eta^{\kappa}_Q = \max\{|\lambda_Q^{\kappa}-\lambda_Q^{\kappa-1}|\},\nonumber\\
	& \eta^{\kappa}_V = \max\{|\lambda_V^{\kappa}-\lambda_V^{\kappa-1}|\}, \quad
	\eta^{\kappa}_I = \max\{|\lambda_I^{\kappa}-\lambda_I^{\kappa-1}|\}.\nonumber
	\end{align}
	\item \emph{Check convergence:} Compare maximum difference with the stopping criteria $\hat\eta$: 
	\begin{equation} 
	\eta_P^{\kappa}\leq\hat{\eta}_P~, \quad
	\eta_Q^{\kappa}\leq\hat{\eta}_Q~, \quad
	\eta_V^{\kappa}\leq\hat{\eta}_V~, \quad
	\eta_I^{\kappa}\leq\hat{\eta}_I~. \label{criterion}
	\end{equation} 
	If the criteria \eqref{criterion} are satisfied: Algorithm converged.\\
	If at least one criterion from \eqref{criterion} is not satisfied: Increase iteration count to $\kappa=\kappa+1$, and move back to step 2).
\end{enumerate}
\vspace{+2pt}
\hrule
\vspace{+6pt}

While the above algorithm is straightforward to implement and offers scalability, it has some drawbacks.

First, since $\alpha$ and $\gamma$ only appear as variables in the uncertainty margins and not as part of the deterministic AC OPF problem, the iterative algorithm requires those parameters to be pre-specified. 
Therefore, the iterative AC CC-OPF algorithm is not able to optimize the system response to uncertainty, defined by $\alpha$ or $\gamma$, in a way that minimizes fluctuations on congested lines or voltage variations at buses with tight voltage constraints. 

Second, the solution is not guaranteed to converge, or to converge to a (locally) optimal solutions.   
In the simulations conducted for this paper, the algorithm converged within a few iterations when $\hat{\eta}_P,~\hat{\eta}_Q,~\hat{\eta}_V$ and $\hat{\eta}_I$ are chosen not smaller than 0.001 MVA for active and reactive power, and $10^{-5}$ p.u. for voltage and current magnitudes. Furthermore, the iterative algorithm found similar solutions as the one-shot optimization algorithm which is guaranteed to converge to a locally optimal solution. Those results indicate that, although the algorithm does not have convergence guarantees, it performs well in practice.



\section{Alternative Chance Constraint Formulations}
\label{sec:Alternative}
When solving the problem using the iterative approach, the uncertainty margins $\lambda_P^{\kappa},~\lambda_Q^{\kappa},~\lambda_V^{\kappa}$ and $\lambda_I^{\kappa}$ are calculated in a step which is separate from solving the AC OPF. Separating the AC OPF solution process from the calculation of the uncertainty margins  enables the use of other methods than the analytical reformulation to compute the uncertainty margins in each iteration. In particular, we are able to apply more accurate, but computationally heavy methods to calculate the margins without sacrificing problem tractability and scalability. Note that although we change the computation of the uncertainty margin, the convergence criterion for the iterative algorithm remains the same, i.e., it terminates when the uncertainty margins no longer change between iterations. 

\subsubsection{Uncertainty Margins from Monte Carlo Simulation}
The analytical uncertainty margins \eqref{currentIII}, \eqref{lambdas} represent an approximate quantile of the current, voltage and generation distributions. By running a Monte-Carlo simulation where we sample the uncertainty vector $\omega$ and calculate the resulting power flows for a large number of samples, we can compute an empirical distribution function and the corresponding quantiles to define Monte Carlo-based uncertainty margins. Since the AC power flow equations are non-linear and the samples might not follow a symmetrical distribution, we calculate the upper and lower uncertainty margins separately.  

For an example voltage constraint, a Monte Carlo simulation is used to determine the distribution function of the voltage around the forecasted solution $v_{i}(\mathbf{x})$. We then determine the upper $(1-\epsilon)$ and lower $(\epsilon)$ quantiles of the distribution, denoted by $v_{i}^{1-\epsilon},~v_{i}^{\epsilon}$, and calculate the constraint tightenings~by 
\begin{equation}
\label{monte_carlo}
\lambda_{V,i}^{U}=v_{i}^{1-\epsilon}-v_{i}(\mathbf{x})\quad \text{and} \quad \lambda_{V,i}^{L}=v_{i}(\mathbf{x})-v_{i}^{\epsilon}.
\end{equation}
Since the outcome of the Monte Carlo simulation and the resulting values of the uncertainty margins depend on the solution to the AC OPF, we need to rerun the Monte Carlo in each iteration $\kappa$. 


\subsubsection{Uncertainty Margins for Joint Chance Constraints}
While the uncertainty margins \eqref{currentIII}, \eqref{lambdas} represent the quantiles of the respective constraints, they can more generally be interpreted as a security margin which is necessary to secure the system against uncertainty. 
A similar constraint tightening can be observed in other stochastic and robust formulations of the OPF problem, e.g. in chance-constrained formulations with joint chance constraints \cite{vrakopoulou2013, vrakopoulou2013AC}, which ensure that all constraints hold \emph{jointly} with a pre-described probability $1-\epsilon_J$. 

To reformulate the joint chance constraints, the above references rely on the scenario approach \cite{campi2006} or a robust version of it \cite{margellos2014}. The scenario approach guarantees that the joint chance constraint will hold if all constraints are satisfied for a defined number $N_S$ of randomly drawn samples. The required number of samples depends on the number of decision variables $N_X$ and the acceptable joint violation probability $\epsilon_J$ \cite{campi2009},
\begin{equation}
\label{scenapp}
N_S \geq \frac{2}{\epsilon_J}\left(\ln\frac{1}{\beta} + N_X \right).
\end{equation}
The value $1-\beta$ represents the confidence level for satisfying the chance constraint, and is typically chosen to be very small. 

Representing uncertainty sets based on samples, 
or representing a robust set through the vertices as in \cite{vrakopoulou2013AC}, 
only applies when the underlying problem is convex, since convexity guarantees that a linear combination of feasible points remain within the convex, feasible set. 
In \cite{vrakopoulou2013AC}, the AC CC-OPF is solved by representing the AC power flow constraints through a convex relaxation. However, while the scenario approach is applicable for the relaxed problem, it does not provide guarantees for the actual AC problem\footnote{The relaxation might not be tight, in which case the solution to the relaxed power flow equations are not physical solutions to the AC power flow problem.}. 

Here, we take a different approach. In each iteration, we first solve \eqref{approxACCCOPF} applying the full, non-convex AC power flow equations
\footnote{We might apply solution methods based on a convex relaxation to find the solution for the AC OPF \eqref{approxACCCOPF}, but would check AC power flow feasibility of the resulting solution. Tightness of the relaxation can be guaranteed, e.g., by applying higher-order moment relaxations as proposed in \cite{molzahn2015}.}. 
%
We use a set of samples $\mathcal{S}$, with $|\mathcal{S}|=N_S$ 
as prescribed by the scenario approach \eqref{scenapp}, and run an AC power flow simulation for each of the scenarios $\hat\omega\in \mathcal{S}$. 
While the overall problem is not convex, we assume that the AC power flow solutions for each of the samples remain in vicinity of the local optimum, where the problem is known to be locally convex and the convexity assumption required for the scenario approach can be expected to hold. 

To ensure that the constraints hold for all samples, we define the uncertainty margin as the difference between the forecasted value of the variable and the highest/lowest observed magnitude, e.g. for the example voltage constraint, 
\begin{subequations}
\label{scenario_approach}
\begin{align}
&\lambda_{V,i}^U=\max_{\hat\omega\in\mathcal{S}}\{v_i(\hat\omega)\}-v_{i}(\mathbf{x}), \\ 
&\lambda_{V,i}^L=v_{i}(\mathbf{x})-\min_{\hat\omega\in\mathcal{S}}\{v_i(\hat\omega)\}.
\end{align}
\end{subequations}

\begin{figure}
	\includegraphics[width=0.95\columnwidth]{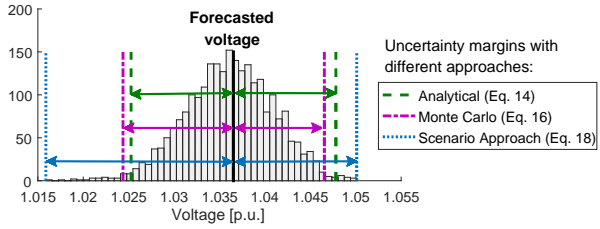}
	\centering
	\caption{Uncertainty margins $\lambda_{V,j}^u$, $\lambda_{V,j}^l$ for an example voltage constraint, as obtained with a given set of samples. The black line represents the forecasted voltage, and the histogram shows the empirical voltage distribution for the given sample set. The uncertainty margins are given by the distance between the black line and the corresponding coloured lines: Analytical reformulation (green), Monte Carlo simulation (purple) and Scenario approach (blue). }
	\label{UncMargDemo} 
\end{figure}


\subsubsection{Comparison of Uncertainty Margins}
Fig. \ref{UncMargDemo} compares the margins calculated with the analytical approach \eqref{currentIII}, \eqref{lambdas}, the Monte Carlo simulation \eqref{monte_carlo} and the scenario approach \eqref{scenario_approach}. The histogram represents the empirical distribution of the voltage magnitudes, as calculated based on the samples used in the Monte Carlo simulation and the scenario approach. 
The analytical margins are symmetric, leading to a larger upper margin and a smaller lower margin compared with the empirical quantiles obtained from the Monte Carlo simulation. 
The scenario approach has much larger uncertainty margins than either of the two other approaches, as these margins are determined by the worst-case among the drawn scenarios. 
The uncertainty margins for the scenario approach are however not directly comparable to the analytical and Monte Carlo approach, since the scenario approach enforces a joint violation probability. In this example, the scenario approach guarantees a joint violation probability of $\epsilon_J\leq0.1$, while the analytical and Monte Carlo approach limit the violation probability of each constraint to $\epsilon\leq0.01$. 


\section{Case Study}
\label{sec:Case}
In the case study, we demonstrate the performance of the AC CC-OPF algorithms and reformulations in terms of operating cost, computational time and chance constraint satisfaction in both in-sample and out-of-sample tests. In addition, we assess the ability of the chance constraints to control the expected size of constraint violation. 

\subsubsection{Investigated approaches}
We compare the following versions of the AC CC-OPF:
\begin{itemize}
    \item[\textbf{A}] \textbf{Standard AC OPF} without consideration of uncertainty.
    \item[\textbf{B}] \textbf{One-Shot AC CC-OPF} with both pre-determined and optimized uncertainty response $\alpha,~\gamma$.
    \begin{itemize}
        \item[\textbf{B1}] \textbf{Fixed} $\boldsymbol{\alpha,\gamma}$ where $\alpha,\gamma$ are pre-defined.
        \item[\textbf{B2}] \textbf{Optimized} $\boldsymbol{\alpha,\gamma}$ where $\alpha$ and $\gamma$ are optimized with the generation dispatch.
    \end{itemize}
    \item[\textbf{C}] \textbf{Iterative AC CC-OPF} with $\alpha,\gamma$ pre-defined, and different uncertainty margin definitions:
    \begin{itemize}
        \item[\textbf{C1}] \textbf{Analytical} with uncertainty margins defined by the closed-form expressions \eqref{lambdas}.
        \item[\textbf{C2}] \textbf{Monte Carlo} with empirical uncertainty margins obtained from a Monte Carlo simulation \eqref{monte_carlo}.
        \item[\textbf{C3}] \textbf{Scenario Approach} with empirical uncertainty margins obtained from the limiting scenarios \eqref{scenapp}.
    \end{itemize}
\end{itemize}

\subsubsection{Test systems}
For all test cases, we assume that uncertainty is observed as fluctuations in the net load, with standard deviations given as a percentage of forecasted load. The uncertainty levels are chosen to obtain congested, but feasible test cases. The uncertainty in the net load represents a combination of uncertainty from sources connected at lower voltage levels, and hence include fluctuations in both load and renewable generation. This type of modelling assumption is made to construct reasonable, reproducible data sets for the stochastic test cases. In general, the framework also enables the modelling of additional uncertain power injections representing, e.g., large wind farms or solar installations connected directly at the transmission level.
For the analytical chance constraint reformulation, we assume that $\omega$ follows a multivariate normal distribution. We do not consider unit commitment, and hence set the lower generation limits to zero. 
We run simulations for four different test systems of different size:

\emph{IEEE RTS96 One Area Test Case} provided with Matpower 5.1 \cite{zimmermann2011}, with the maximum generation limits increased by a factor of 1.5. All 17 loads are uncertain, with standard deviations equal to 10\% and zero correlation between loads. 

\emph{IEEE 118 Bus Test System} from the NICTA Energy System Test Archive \cite{coffrin2014}. The generation limits are increased 1.5, and the system is split into three zones as in \cite{roald2016corrective}. We assume that all 99 loads are uncertain with standard deviation of 5\%, a correlation coefficient of $\rho=0.3$ within each zone, and zero correlation between loads in different zones. 

\emph{IEEE 300 Bus Test System} from the NICTA Energy System Test Archive \cite{coffrin2014}. Loads with consumption between 0 and 100 MW are assumed to be uncertain, with standard deviations of 5\% and zero correlation. This corresponds to 131 uncertain loads, and 21\% of the total system demand.  

\emph{Polish 2383 Bus, Winter Peak Test Case} from Matpower 5.1 \cite{zimmermann2011}. The upper generation limits for active power are increased by a factor of 2, and the reactive power capability of generators at PV buses is increased by +/- 10 MVAr for all generators.
All loads with consumption between 10 and 50~MW are assumed to be uncertain, with standard deviations equal to 10\% and zero correlation. This corresponds to 941 uncertain loads, and 67\% of the total system demand.  

As a base case, we enforce all chance constraints with violation probabilities $\epsilon_P = \epsilon_Q = \epsilon_V = \epsilon_I = 0.01$.
For the iterative solution algorithms we set $\hat{\eta}_P=\hat{\eta}_Q=0.001 MVA$, $\hat{\eta}_V=10^{-5} p.u.$ and $\hat{\eta}_I=0.001 kA$.
When $\alpha,\gamma$ are pre-defined, we assume that $\gamma$ is the power ratio of the forecasted load, and that the participation factors $\alpha$ are given by  $\alpha_i=p_{G,i}^{max}/\sum_{g\in\mathcal{G}} p_{G,g}^{max}$ as in \cite{roald2013}.

\subsection{Comparison of Solution Algorithms}
We first compare the solutions obtained with the deterministic AC OPF (A), the iterative approach (C1), the one-shot optimization with fixed participation factors $\alpha,\gamma$ (B1) and the one-shot optimization with $\alpha,\gamma$ as decision variables (B2) in terms of cost and solution time. 
The computational times are based on our own AC OPF implementation, 
solved using KNITRO. 
The results are listed in Table \ref{table_algorithms}. 

{\begin{table} 
\caption{Cost of generation (relative to the deterministic OPF) and solution time with different solution approaches.} 
\label{table_algorithms}
\centering
\begin{tabular}{lllll}
                & \textbf{Standard} & \textbf{Iterative} & \textbf{O-S Fixed}    & \textbf{O-S Opt}\\
                & \textbf{(A)}           & \textbf{(C1)}      & \textbf{(B1)}              & \textbf{(B2)} \\[+4pt]

\textbf{RTS96}           & & & & \\[+1pt]
~Cost            & \$36'771          & +7.7\%        & +7.6\%            & +3.2\%\\[+1pt]
~Time            & 1.5s              & 4.3s          & 4.3s              & 12.5s \\[+1 pt]
\textbf{118 Bus}         & & & & \\[+1pt]
~Cost            & \$3'504.1         & +1.9\%        & +1.9\%            & -\\[+1pt]
~Time            & 2min 5s           & 7min 46s      & 10min 51s         & - \\[+1pt]
\end{tabular}
\end{table}}

Due to the introduction of the uncertainty margin, the AC CC-OPF formulations have higher operational cost than the deterministic solution. 
We observe that the iterative approach (C1) and the one-shot optimization with fixed response $\alpha,~\gamma$ lead to similar solutions, while the one-shot solution with optimized $\alpha, \gamma$ has lower cost. Co-optimizing $\alpha, \gamma$ does however significantly increase computational complexity and solution time. Already for the modestly sized 118 bus system, we observe that the one-shot optimization (both B1 and B2) have longer solution times than the iterative approach. The optimization (B2) does not converge at all, even after running for 13h on a desktop computer.

Another important observation from this comparison is that the solution obtained with the iterative (C1) and one-shot (B1) solution approaches are very close to each other. They do not only have a similar cost, as shown in Table \ref{table_algorithms}, but also finds a similar optimal point.


\begin{table} 
\caption{For the four different test cases (listed with system size and number of uncertain loads): Solution times, number of iterations and generation costs in the first, second and final iterations.}
\label{table_solutions}
\centering
\begin{tabular}{lllll}
                & \textbf{RTS96} & \textbf{118 Bus} & \textbf{300 Bus}  & \textbf{Polish}\\[+6pt]
Buses            & 24          & 118            & 300             & 2383\\
Uncertain loads  & 17          & 99             & 131             & 941\\[+4pt]
Solution time            & 0.54s          & 1.15s            & 3.37s             & 31.89s\\[+4pt]
Iterations      & 5              & 4                & 5                 & 4 \\[+4pt]
Cost (1st)      & 36 771         & 3504.1           & 16 779            & 787 987\\
Cost (2nd)      & 40 249         & 3575.7           & 17 173            & 802 529\\
Cost (final)    & 40 127         & 3575.3           & 17 143            & 802 238\\
\end{tabular}
\vspace{-8pt}
\end{table}

\subsection{Scalability of the Iterative Approach}
We now demonstrate how the iterative approach addresses the issue of scalability by utilizing existing OPF tools. 
We implement the iterative AC CC-OPF using the standard Matpower 5.1 ``runopf" function \cite{zimmermann2011} with the default MIPS solver to solve the deterministic OPF, and run the problem for all four test systems. The resulting times and the number of iterations are shown in Table \ref{table_solutions}. We also show the evolution of the generation cost between first, second and last iterations.


The problems converge within 4-5 iterations, and the solutions are obtained within half a minute on a standard desktop computer, even for the Polish test case with 941 uncertain loads and 2383 buses. 
Further, the main change in cost happens between the first and the second iteration, with only minor adjustments until final convergence. This implies that even if the problem does not converge, the solution at intermediate iterations might already perform well.

\subsection{Evaluation of Chance Constraint Reformulation}
The AC CC-OPF with the analytical reformulation has two main sources of inaccuracy. 
First, the impact of uncertainty is approximated by a linearization. 
Second, the analytical reformulation assumes a normal distribution, which might be an inaccurate description of the true uncertainty distribution. 
To assess the accuracy of the analytical reformulation, we solve the problem and compare the pre-described acceptable violation probability with observed violation probabilities from a Monte Carlo simulation. The comparisons are based on the RTS96 system, and the AC CC-OPF is solved using the iterative approach with Matpower and MIPS.


\subsubsection{Accuracy of the Linearization: In-Sample Test}
To assess the accuracy of the linearization, we assume perfect knowledge of the distribution and run an in-sample test based on 10'000 samples from the assumed multivariate normal distribution. 
Based on these samples, we assess the empirical violation probabilities through a Monte Carlo simulation. We perform the assessment for different standard deviations $\sigma_W=\{0.075, 0.1, 0.125\}$ corresponding to fluctuations of different size, and different acceptable violation probabilities $\epsilon=\epsilon_V=\epsilon_I=\{0.01, 0.05, 0.1\}$ where a smaller $\epsilon$ implies the estimation of quantiles that are further into the tail of the distribution and thus further away from the operating point. 
The acceptable violation probabilities of the generator constraints on active and reactive power are kept constant at $\epsilon_P=\epsilon_Q=0.01$, as these constraints are not heavily influenced by the linearization of the power flow equations. 

Table \ref{insample} shows the results of the in-sample testing, with the maximum observed empirical violation probability $\epsilon_{emp}$ for any constraint, as well as the joint violation probability (the percentage of samples that have at least one violated constraint). The upper part of the table shows the results for different standard deviations, and the lower part shows the results for varying values of $\epsilon$. Since the linearization is a better approximation close to the forecasted operating point, we observe that the maximum empirical violation probability is closer to the acceptable value for small standard deviations $\sigma_W$ and large acceptable violation probabilities $\epsilon$. For large standard deviations and small $\epsilon$, the maximum violation probability is higher than the acceptable value ($\epsilon_{emp}>\epsilon$) implying a violation of the chance constraint, whereas the method actually lead to conservative solutions ($\epsilon_{emp}<\epsilon$) for larger values of $\epsilon$. However, the empirical violation probability $\epsilon_{emp}$ is always within $\pm 0.01$ of the pre-described acceptable $\epsilon$.


We further observe that the AC CC-OPF limits the joint violation probability $\epsilon_J$ (i.e., the probability of observing at least one violation) to a relatively small percentage, even though it only aims at enforcing the separate chance constraints. 

{\begin{table}
\caption{In-sample testing of the analytical chance constraint (C1): Maximum observed violation probability $\epsilon_{emp}$, and observed joint violation probability $\epsilon_J$ for different standard devations $\sigma_W$ (with $\epsilon=0.01$) and acceptable violation probabilities $\epsilon$ (with $\sigma_W = 0.1$).}
\label{insample}
\centering
\begin{tabular}{lccc}
                        & $\sigma_W=0.075$ & $\sigma_W=0.1$ & $\sigma_W=0.125$ \\[+4pt]
Max. $\epsilon_{emp}$   & 0.011  & 0.013        & 0.017                \\[+3pt]
Joint $\epsilon_J$      & 0.065  & 0.065        & 0.081      \\
& & & \\[+4pt]
                        & $\epsilon=0.01$ & $\epsilon=0.05$ & $\epsilon=0.1$ \\[+4pt]
Max. $\epsilon_{emp}$   & 0.013           & 0.044           & 0.092                \\[+3pt]
Joint $\epsilon_J$      & 0.065           & 0.137        & 0.219      \\
                   
\end{tabular}
\end{table}}

\subsubsection{Accuracy of Normal Distribution: Out-Of-Sample Test}

To check whether the assumption of a normal distrbution leads to accurate results, we perform an out-of-sample test with power injection samples based on the historical data from Austrian Power Grid. We assign one set of historical samples, in total 8492 data points,  to each uncertain load. The samples are then rescaled to match the assumed standard deviation $\sigma_W$ and assumed correlation coefficient $\rho=0$. We use 5000 samples for the evaluation of constraint violation probabilities.

The results of the out-of-sample testing based on APG historical data is shown in Table \ref{outofsample}, including both the maximum violation probability for any individual constraint and the joint violation probability as above. Since the out-of-sample test includes inaccuracies both due to linearization errors and non-normally distributed samples, the maximum observed violation probability is higher in the out-of-sample test than in the in-sample test. The difference is however not particularly large, and the empirical violation probability $\epsilon_{emp}$ is still within $\pm 0.01$ of the pre-described acceptable $\epsilon$. Since the inaccuracy due to the distribution assumption is on the same order of magnitude as the inaccuracy due to the linearization error, the normal distribution appears to be a reasonable model for current and voltage magnitudes, even though the power injections are not normally distributed. 



Similarly, the joint violation probability observed with the APG data is slightly higher than in the in-sample test, but remains in the same range.

{\begin{table} 
\caption{Out-of-sample testing of the analytical chance constraint (C1): Maximum observed violation probability $\epsilon_{emp}$, and observed joint violation probability $\epsilon_J$ for different standard devations $\sigma_W$ (with $\epsilon=0.01$) and acceptable violation probabilities $\epsilon$ (with $\sigma_W = 0.1$).}
\label{outofsample}
\centering
\begin{tabular}{lccc}
                    & $\sigma_W=0.075$ & $\sigma_W=0.1$ & $\sigma_W=0.125$ \\[+4pt]
Max. $\epsilon_{emp}$ & 0.017          & 0.014        & 0.020  \\[+2pt]
Joint $\epsilon_J$      & 0.075  & 0.074        & 0.092      \\
& & & \\[+4pt]
                    & $\epsilon=0.01$ & $\epsilon=0.05$ & $\epsilon=0.1$ \\[+4pt]
Max. $\epsilon_{emp}$ & 0.014           & 0.054           & 0.093 \\[+2pt]
Joint $\epsilon_J$      & 0.074  & 0.145        & 0.233      \\
\end{tabular}
\vspace{-6pt}
\end{table}}

\subsection{Controlling Violation Size Through Chance Constraints}
\label{sec:ExpectedViol}
The AC CC-OPF only controls the probability of constraint violations, without further consideration of their size or impact on system operation. However, it is possible to assess the impact of the acceptable violation probability $\epsilon=\epsilon_P =\epsilon_Q =\epsilon_I =\epsilon_V$ through an a-posteriori simulation. In this example, we assess how $\epsilon$ influences the size of the constraint violations. To compute the size of the violations, we solve the AC CC-OPF using the analytical chance constraint formulation with $0.01 \leq \epsilon \leq 0.15$ in steps of $\delta \epsilon = 0.01$. We then calculate the expected size of the constraint violations for each constraint, based on a Monte Carlo simulation with 2000 samples drawn from the multivariate normal distribution. The size of the constraint violation is defined as the amount by which the limit is exceeded, and is set to zero for samples where the constraints are not violated. 

\begin{figure}
	\includegraphics[width=0.97\columnwidth]{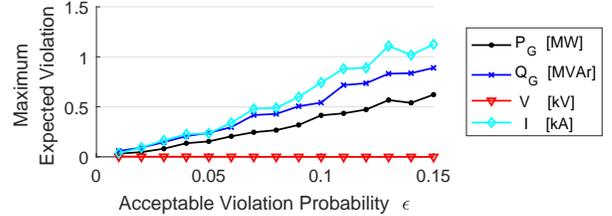}
	\centering
	\caption{Empirically observed size of violation for a range of different acceptable violation probabilities $\epsilon$, as calculated through a Monte Carlo simulation with 2000 multivariate normal samples. The figure shows the maximum observed expected violation probability across all constraints, considering four constraint categories: Generator active power constraints $P_G$, generator reactive power constraints $Q_G$, voltage magnitude constraints $V$ and current magnitude constraints $I$.}
	\label{ViolationSize} 
\end{figure}

Fig. \ref{ViolationSize} shows maximum expected violation size observed among the constraints, for the different choices of $\epsilon$. Since the different types of constraints have different units, we separate the results for the generator active and reactive power constraints, as well as the voltage\footnote{Note that the maximum expected voltage magnitude violation is close to zero in all cases, as there are no active voltage constraints at PQ buses in this case. Therefore, we observe only very few and very small voltage magnitude violations, leading to expecteations close to zero.} and current magnitudes. Not surprisingly, the maximum expected violation increase as the acceptable violation probability increases. However, the increase is not monotonic, indicating that there is not a one-to-one relation between the acceptable violation probability and the resulting violation size. 
Therefore, while reducing the acceptable violation probability generally leads to smaller violations, it is more effective to control the size of the constraint violation directly using, e.g., weighted chance constraints \cite{roald2015}.

\subsection{Comparison with Sample-Based Reformulations}
Finally, we compare the analytical reformulation with the sample-based methods to define the uncertainty margins. 
For the analytical approach (C1), we assume a normal distribution and enforce $\epsilon=0.01$. For the Monte Carlo simulation (C2), we enforce $\epsilon=0.01$ and use 1000 samples from the APG data to calculate the empirical uncertainty margins. For the scenario approach (C3), we enforce $\epsilon_J=0.1$. This corresponds to a prescribed number of $N_S=2465$ samples, which were taken from the APG data. 
To evaluate the actual violation probability, we use 5000 samples from the APG historical data.


In Table \ref{table_joint}, the generation cost, solution time and number of iterations are listed. We also show the maximum empirical violation probability $\epsilon_{emp}$ and joint violation probability $\epsilon_J$. 

{\begin{table} 
\caption{Comparison of the iterative solution approach uncertainty margins based on (C1) the analytical reformulation and $\epsilon=0.01$, (C2) a Monte Carlo simulation and $\epsilon=0.01$, and (C3) the scenario approach and an acceptable joint violation probability of $\epsilon_J=0.1$.}
\label{table_joint}
\centering
\begin{tabular}{llll}
                    & \textbf{Analytical} & \textbf{Monte Carlo}  & \textbf{Scen. Approach}\\
                    & \textbf{(B1)}           & \textbf{(B2)}             & \textbf{(B3)} \\[+6pt]
Cost                & 40 127         & 39 922           & 41 455\\[+4pt]
Max. $\epsilon_{emp}$ & 0.014  & 0.019            & 0.003  \\
Joint $\epsilon_J$  & 0.074            & 0.10              & 0.007             \\[+4pt]
Time                & 0.54s          & 2min 7s          & 6min 12s\\
Iterations          & 5              & 5                & 6 \\
\end{tabular}
\end{table}}

We observe that the analytical (C1) and Monte Carlo (C2) approaches lead to relatively similar solutions. The Monte Carlo solution has slightly lower cost, but higher violation probabilities. The Monte Carlo thus leads to a more significant violation of the acceptable violation probability than the analytical reformulation, despite being based on the full AC power flow equations and requiring no explicit assumption about the distribution. This highlights the sensitivity of the Monte Carlo approach to the availability sufficiently many, high quality samples. A larger number of samples would lead to a better solution, but would also increase computational time. The Monte Carlo already requires a lot more computational resources than the analytical reformulation, with a solution time of 2min 7s compared with 0.5s.


The scenario approach (C3) has a higher cost than the other two solutions, but also significantly lower violation probabilities. 
While the pre-described $\epsilon_J=0.1$, the actual joint violation probability was only $\epsilon_{J,emp}=0.007$. This shows an important drawbacks of the scenario approach: While it guarantees chance constraint feasibility, the solution might be very conservative and possibly far from cost optimal.


\section{Discussion of Practical Implementation Aspects}
\label{sec:Practical}
In the following, we comment on some practical aspects relevant to implementation in industrial applications.
\subsection{Including N-1 Constraints}
The AC CC-OPF model in this paper only considers the system in its pre-contingency state, i.e., with all components in operation. In practice, power systems operation typically requires N-1 security, which guarantees satisfaction of all constraints despite outage of any single component. Extending the above problem to include N-1 security constraints by formulating an AC CC-Security Constrained OPF (AC CC-SCOPF) is possible without major changes to the framework, but requires some conceptual considerations and can be computationally challenging. 

Conceptually, the post-contingency N-1 constraints can be enforced in a similar way as the pre-contingency constraints, by introducing an additional set of AC power flow constraints for each post-contingency state and the corresponding chance constraints on generation, voltages and currents. 
However, an important modelling aspect is the choice of acceptable violation probability for the post-contingency constraints, which might be chosen separately for each contingency\footnote{Indeed, all violation probabilities can be chosen independently, although we have chosen to assume similar violation probabilities $\epsilon_P,~\epsilon_Q, ~\epsilon_V$ and $\epsilon_I$ for each type of constraint}. In literature, \cite{roald2016corrective} assume a similar violation probability for both pre- and post-contingency chance constraints, while \cite{vrakopoulou2013AC} guarantees that both pre- and post-contingency constraints will be satisfied with a pre-defined, joint violation probability.
In systems where the regulator requires N-1 security at all times, this choice of acceptable post-contingency violation probability is reasonable, and is closely linked to our interpretation of chance constraint violations as the probability that the operator needs to take additional control actions in real time. 

A more sophisticated method for choosing acceptable post-contingency violation probabilities is provided in \cite{hamon2013}. Here, the acceptable violation probability for each contingency case is related to an overall probability of system failure, and the acceptable post-contingency violation probability is weighted by the probability of the N-1 contingency. While such an approach is theoretically appealing, there is a chance that the probability of the contingency is so small that the optimization problem allows for very high probability of constraint violation for the case when contingency occurs. Since high probability of violation is typically linked to large expected violations, as demonstrated in Section \ref{sec:ExpectedViol}, this can have potentially catastrophic impact on the system.
To accurately model the risk related to contingencies, it is therefore not sufficient to only consider probability of constraint violation, but also the resulting impact of the violation. If such tools are not available, it might therefore be better to use simple rules of thumbs as in \cite{sundar2016}, where the acceptable violation probability was increased by a factor of two for the post-contingency constraints.

Computationally, N-1 security constraints pose a challenge due to the large number of considered system states, which comes in addition to the consideration of uncertainty. Ref. \cite{roald2016corrective} developed an algorithm based on constraint generation which handles large scale systems with security and chance-constraints, and although the results were limited to the DC power flow case similar ideas based on constraint generation could be applied to the AC case. Related ideas to handle uncertainty and security constraints were also proposed in, e.g., \cite{capitanescu2012}, where the set of considered uncertainty scenarios was discovered sequentially. Furthermore, the proposed iterative scheme provides an interesting path to solve AC CC-SCOPF for large and practical systems, as it can utilize existing implementations to solve the AC SCOPF in combination with an appropriate tool for calculation of the uncertainty margins.

\subsection{Compatibility with Existing Tools}
The iterative solution algorithm for the AC CC-OPF offers significant potential for practical implementation. On the one hand, practical implementation is enabled by the ability to implement the iterative algorithm with existing, industrial AC SCOPF implementations as a core component, without requiring major changes to the solution approach. 
On the other hand, the iterative approach enables sample-based assessment of the uncertainty impact, which can be performed based on power flow solvers already in use in practical operation. This enables the treatment of complicated non-smooth processes such as consideration of PQ-PV bus switching.

An example of the amenability of the iterative approach to practical applications is available from the UMBRELLA project \cite{morales2016}. In this project, an AC SCOPF which used Monte Carlo-based uncertainty margins to account for uncertainty was implemented as part of the developed prototype toolbox. The implemented approach calculated the uncertainty margins based only on an estimated operating point, and did not include any subsequent updates or iterations. Therefore, the uncertainty margins did not provide any guarantees for the final solution. However, despite this  inaccuracy, the approach was shown to significantly reduce the probability of constraint violations. 

\subsection{Availability of Data}
One main problem facing practical implementation of chance-constrained AC OPF (or any other stochastic OPF method) is the availability of high quality data or probabilistic forecasts describing the uncertainty, which in particular must capture the geographical correlation of different uncertainty sources. Data availability is not only an academic problem, but also a practical challenge, as it requires gathering, storing and processing large amounts of information. The challenge of collecting this data was recognized in one of the common recommendations to the ENTSO-E (the European network of transmission system operations for electricity) from the UMBRELLA and iTesla research projects \cite{umbrellaD6p3}.

\section{Conclusions}
\label{sec:Conclusion}
In this paper, we describe a model for the AC CC-OPF problem, and discuss different reformulations and solution algorithms to obtain a tractable optimization problem formulation. 
We first discuss an analytical method to reformulate the chance constraints into closed-form expressions. The reformulation method uses the full, non-linear AC power flow equations for the forecasted operating point, and a linearization around this point to model the impact of uncertainty. 
This partial linearization leads to a linear dependence on the random fluctuations and enables an analytical chance constraint reformulation. We show that this analytical reformulation approach is reasonably accurate for the AC CC-OPF, and that reformulating the constraints assuming normally distributed voltages and currents can provide a good approximation for even cases where the uncertain power injections are not normally distributed.

We discuss different solution algorithms for the problem. The analytically reformulated AC CC-OPF admits direct solution in a one-shot optimization, which allows us to co-optimize the reaction to the fluctuations and hence reduce cost. However, the computational complexity of the one-shot solution easily becomes prohibitive. 
To increase computational tractability for large problems, we discuss an iterative solution method. The iterative algorithm is based on the observation that the uncertainty leads to a constraint tightening, the so-called uncertainty margins, and alternates between (i) solving a deterministic AC OPF problem with fixed uncertainty margins, and (ii) computing the ``correct" uncertainty margins for the resulting operating point. Separating those two steps has two main advantages:

First, it can utilize existing solvers for deterministic AC OPF to enable scalable implementations, as demonstrated in the case study, where we solve the AC CC-OPF for the Polish test case with 941 uncertain loads in 30s on a desktop computer. 

Second, it enables the implementation of more general chance-constraint reformulations in the uncertainty assessment step without sacrificing computational tractability. 

Utilizing the latter characteristic, we propose two alternative, sample-based reformulations of the chance constraints, based on Monte Carlo simulations and the scenario approach. These reformulations require neither approximation or relaxation of the AC power flow equations, nor any limiting assumptions about the distribution beyond the availability of samples. We used the alternative reformulations to benchmark the performance of the analytical reformulation, and showed that the analytical reformulation had lower solution times. Moreover, the analytical reformulation outperformed the Monte Carlo reformulation in terms of enforcing the chance constraint, highlighting the sensitivity of sampling based solutions to the choice of samples. The scenario approach provided the most rigorous guarantees for chance constraint satisfaction, but at the expense of very high cost.

There are several directions for future work in which the proposed method can be improved. First, demonstrating the practicability of the algorithm by including N-1 constraints (as discussed in Section \ref{sec:Practical}) or by extending the modelling to account for the cost of balancing energy is one important direction. Second, the iterative algorithm opens for extensions to formulations beyond chance-constraints, such as weighted chance constraints \cite{roald2015} that limit the size of violations, or the application to robust AC OPF. Third, while the iterative algorithm has been observed to behave well on standard test cases, additional work to develop extensions that provide guarantees for convergence and solution quality are required.

\bibliographystyle{IEEEtran}
\bibliography{Diss_bib}



\end{document}